\pdfoutput=1
\RequirePackage{ifpdf}
\ifpdf 
\documentclass[pdftex]{sigma}
\else
\documentclass{sigma}
\fi

\numberwithin{equation}{section}

\newtheorem{Theorem}{Theorem}[section]
{\theoremstyle{definition}
\newtheorem{Example}[Theorem]{Example}
\newtheorem{Remark}[Theorem]{Remark}
}

\begin{document}

\allowdisplaybreaks

\renewcommand{\PaperNumber}{041}

\FirstPageHeading

\ShortArticleName{On the Linearization of Second-Order ODEs to the Laguerre Form}

\ArticleName{On the Linearization of Second-Order Ordinary\\
Dif\/ferential Equations to the Laguerre Form\\ via Generalized
Sundman Transformations}

\Author{M.~Tahir MUSTAFA, Ahmad Y.~AL-DWEIK and Raed A.~MARA'BEH}
\AuthorNameForHeading{M.T.~Mustafa, A.Y.~Al-Dweik and R.A.~Mara'beh}

\Address{Department of Mathematics \& Statistics, King Fahd
University of Petroleum and Minerals,\\ Dhahran 31261, Saudi Arabia}
\Email{\href{mailto:tmustafa@kfupm.edu.sa}{tmustafa@kfupm.edu.sa},
\href{mailto:aydweik@kfupm.edu.sa}{aydweik@kfupm.edu.sa}, \href{mailto:raedmaraabeh@kfupm.edu.sa}{raedmaraabeh@kfupm.edu.sa}}
\URLaddress{\url{http://faculty.kfupm.edu.sa/math/tmustafa/},\\
\hspace*{10.5mm}\url{http://faculty.kfupm.edu.sa/MATH/aydweik/}}

\ArticleDates{Received February 16, 2013, in f\/inal form May 25, 2013; Published online May 31, 2013}

\Abstract{The linearization problem for nonlinear second-order
ODEs to the Laguerre form by means of generalized Sundman
transformations (S-transformations) is considered, which has been
investigated by Duarte et al.\ earlier. A characterization of these
S-linearizable equations in terms of f\/irst integral and procedure
for construction of linearizing S-transformations has been given
recently by Muriel and Romero. Here we give a new characterization of
S-linearizable equations in terms of the coef\/f\/icients of ODE and
one auxiliary function. This new criterion is used to obtain the
general solutions for the f\/irst integral explicitly, providing a
direct alternative procedure for constructing the f\/irst integrals
and Sundman transformations. The ef\/fectiveness of this approach is
demonstrated by applying it to f\/ind the general solution for
geodesics on surfaces of revolution  of constant curvature in a
unif\/ied manner.}


\Classification{34A05; 34A25}

\section{Introduction}

The mathematical modeling of many physical phenomena leads to such
nonlinear ordinary dif\/ferential equations (ODEs) whose analytical
solutions are hard to f\/ind directly. Therefore, the approach of
investigating nonlinear ODEs via transforming to simpler ODEs
becomes important and has been quite fruitful in analysis of
physical problems. This includes the classical linearization
problem of f\/inding transformations that linearize a given ODE. For
the linearization problem of second-order ODEs via point
transformations, it is known that these must be at most cubic in
the f\/irst-order derivative and its coef\/f\/icients should satisfy the
Lie linearization test
\cite{Ibragimov2006,Ibragimov2004,fazal,Meleshko}. The
implementation of the Lie linearization method requires solving
systems of partial dif\/ferential equations (PDEs). It is also well
known that only second-order ODEs admitting 8 dimensional Lie
symmetry algebra
 pass the Lie linearization test, which makes it a~restricted class of ODEs.
 In order to consider a larger class of ODEs, linearization
 problem via nonlocal transformations has been investigated in~\cite{Chandrasekar,Duarte,Euler}. Many of these transformations are of the form
\begin{gather}\label{OP}
u(t) =\psi  ( x,y   ), \qquad    dt=\phi  ( x,y,y'  ) dx, \qquad \psi_y \phi\neq0,
\end{gather}
and the linearization problem via transformations (\ref{OP}), in
general,  is an open problem. In case that $\phi=\phi  ( x,y
 )$, the transformations of type~(\ref{OP}) are called
generalized Sundman transformations~\cite{Euler2003} and equations
that can be linearized by means of generalized Sundman
transformations to the Laguerre form $u_{tt}=0$ are called
S-linearizable~\cite{MR2010a}. These transformations have also
been utilized to def\/ine Sundman symmetries of ODEs~\cite{Euler1997, Euler,Euler2003}. It should be mentioned that
another special classes of nonlocal transformations of type~(\ref{OP}) with polynomials of f\/irst degree in~$y'$ for~$\phi
 ( x,y,y'  )$ have also been studied  in~\cite{Chandrasekar2006,MR2011}.

Duarte et al.~\cite{Duarte} showed that the S-linearizable
second-order equations
\begin{gather}\label{a}
y''=f(x,y,y')
\end{gather}
are at most quadratic in the f\/irst derivative, i.e.\ belong to the
family of equations of the form
\begin{gather}\label{a0}
y''+F_2(x,y){y'}^2+F_1(x,y)y'+F(x,y)=0.
\end{gather}
Precisely, the free particle equation
\begin{gather*}
 u_{tt}=0
\end{gather*}
can be transformed by an arbitrary generalized Sundman
transformation
\begin{gather}\label{a2}
u(t) =\psi \left( x,y  \right), \qquad dt=\phi \left( x,y  \right) dx, \qquad \psi_y \phi\neq0,
\end{gather}
to the family of equations of the form (\ref{a0}) with the
coef\/f\/icients $F(x,y)$, $F_1(x,y)$ and $F_2(x,y)$ satisfying the
following system of partial dif\/ferential equations
\begin{gather}\label{a3}
A{F_2} = {A_y},\qquad A{F_1} = {B_y} + {A_x},\qquad  AF = {B_x},
\end{gather}
where
\begin{gather}\label{a4}
A = \frac{{{\psi _y}}}{\varphi }, \qquad B =\frac{{{\psi _x}}}{\varphi }.
\end{gather}
 They also gave a characterization of these
S-linearizable equations in terms of the coef\/f\/icients. Muriel and
Romero \cite{MR2010a} further studied S-linearizable equations and
proved that these must admit f\/irst integrals that are polynomials
of f\/irst degree in the f\/irst-order derivative.

\begin{Theorem}[\cite{MR2010a}] \label{t1}
The ODE \eqref{a} is S-linearizable if and only if it admits a
first integral of the form $w(x,y,y')=A(x,y)y'+B(x,y)$. In this
case ODE has the form \eqref{a0}. If a linearizing
S-transformation \eqref{a2} is known then a first integral
$w(x,y,y')=A(x,y)y'+B(x,y)$ of \eqref{a0} is defined by~\eqref{a4}. Conversely, if a first integral
$w(x,y,y')=A(x,y)y'+B(x,y)$ of \eqref{a0} is known then a
linearizing S-transformation can be determined by
\begin{gather*}
\psi(x,y) =\eta(I(x,y)) ,\\
\phi(x,y) =\frac{{\psi}_y}{A} \qquad \text{or} \qquad \phi(x,y) =\frac{{\psi}_x}{B} \quad \text{if} \quad B\neq 0,
\end{gather*}
where $I(x,y)$ is the first integral of
\begin{gather}\label{a42}
y' =-\frac{B}{A}.
\end{gather}
\end{Theorem}

Moreover, Muriel and Romero in \cite{MR2010a,MR2009a} revisited
Duarte results \cite{Duarte}, presented the following equivalent
characterization of S-linearizable ODE of the form (\ref{a0}), and
 also provided constructive methods, as given in Theorem \ref{th-1}, to derive the linearizing
S-transformations.
\begin{Theorem}[\cite{MR2010a}] \label{t2}
Let us consider an equation of the form \eqref{a0} and let ${S_1}$
and ${S_2}$ be the functions defined by
\begin{gather}
{S_1}(x,y) = {{F_1}_{y}} - 2{{F_2}_{x}},\nonumber\\
{S_2}(x,y) = {({F}{F_2} + {F_{y}})_y} + {({{F_{2}}_x} -{{F_{1}}_y})_x} + ({{F_{2}}_x} - {{F_{1}}_y}){F_1}. \label{a5}
\end{gather}
The following alternatives hold:
\begin{itemize}\itemsep=0pt
\item If ${S_1}=0$ then equation \eqref{a0} is S-linearizable
if and only if  ${S_2}=0$

\item If ${S_1} \ne 0$, let  ${S_3}$ and ${S_4}$ be the functions
defined by
\begin{gather}
{S_3}(x,y) = {\left( {\frac{{{S_2}}}{{{S_1}}}} \right)_y} -({{F_{2}}_x} - {{F_{1}}_y}),\nonumber\\
{S_4}(x,y) = {\left( {\frac{{{S_2}}}{{{S_1}}}} \right)_x} +{\left( {\frac{{{S_2}}}{{{S_1}}}} \right)^2} + {F_1}\left({\frac{{{S_2}}}{{{S_1}}}} \right) + {F}{F_2} + {F_{y}}. \label{a6}
\end{gather}
Equation \eqref{a0} is S-linearizable if and only if ${S_3}=0$ and
${S_4}=0.$
\end{itemize}
\end{Theorem}

\begin{Theorem}[\cite{MR2010a}] \label{th-1}
Consider an equation of the form \eqref{a0} and let $S_{1}$ and
$S_{2}$ be the functions defined by \eqref{a5}. The following
alternative hold:
\begin{itemize}\itemsep=0pt
\item If $S_{1}=0$ then the equation has a first integral of the
form $w = A(x,y)y' + B(x,y)$ if and only if $S_{2}=0$. In this
case~$A$ and~$B$ can be given as $A = q  {\rm e}^P$, $B = Q $, where
$P$ is a~solution of the system
\begin{gather}\label{eq-20}
 {P_x} = \frac{1}{2}{F_1},\qquad {P_y}
= {F_2},
\end{gather}
$q$ is a nonzero solution of
\begin{gather}\label{eq-21}
q''(x) + f(x)q(x) = 0,
\end{gather} where
\begin{gather*}
f(x) = F{F_2} + {F_y} - \frac{1}{2}{F_{1x}} - \frac{1}{4}{F_1}^2
\end{gather*}
and $Q$ is a solution of the system
\begin{gather*}
{Q_x} = F q {\rm e}^P  ,\qquad {Q_y}  = \left({\frac{1}{2}{F_1} - \frac{{q'}}{q}} \right)q{\rm e} ^P.
\end{gather*}
\item If ${S_1} \ne 0$ then the equation has a first integral has
a first integral of the form $w = A(x,y)y' + B(x,y)$ if and only
if ${S_3} = {S_4} = 0$, where ${S_3}$ and  ${S_4}$ are the
functions defined by~\eqref{a6}. In this case $A$ and~$B$ can be
given as $A = {\rm e}^P $, $B = Q $, where $P$ is a solution of the
system
\begin{gather}\label{eq-25}
{P_x} = {F_1} + \frac{{{S_2}}}{{{S_1}}} ,\qquad {P_y} =
{F_2},
\end{gather}
and $Q$ is a solution of the system
\begin{gather*}
{Q_x} = F {\rm e}^P  ,\qquad {Q_y} =  -
\left( {\frac{{{S_2}}}{{{S_1}}}} \right) {\rm e}^P .
\end{gather*}
\end{itemize}
\end{Theorem}

In this paper, a new characterization of S-linearizable equations
in terms of the coef\/f\/icients and one auxiliary function is given,
and the  equivalence with the old criteria is proved. This
criterion is used to provide explicit general solutions for the
auxiliary functions~$A$ and~$B$ given in~(\ref{a4}) which can be
directly utilized to obtain the f\/irst integral of~(\ref{a0}). So,
using Theorem~\ref{t1}, the linearizing generalized Sundman
transformations can be constructed by solving the f\/irst-order ODE~(\ref{a42}). The method is illustrated in examples where we
recover the Sundman transformations of Muriel and Romero in~\cite{MR2010a}.

As an application, we express the system of geodesic equations for
surfaces of revolution  as a single second-order ODE and use our
method to f\/ind the general solution for geodesics on surfaces of
revolution  of constant curvature in a unif\/ied manner.

In this paper, we have focused on S-linearization to the Laguerre
form $u_{tt}=0$. For an account of S-linearization  to any linear
second-order ODE, the reader is referred to~\cite{Warisa}.

\section{The method for constructing the f\/irst integrals\\ and Sundman transformations}

When the ODE (\ref{a0}) is S-linearizable, Theorem~\ref{t2} does
not give a method to construct the linearizing generalized Sundman
transformations. In order to derive a method to obtain linea\-ri\-zing
generalized Sundman transformations~(\ref{a2}) of a given
S-linearizable equation (\ref{a0}), Muriel and Romero
\cite{MR2010a} found additional relationships between the
functions $\phi$ and $\psi$ in~(\ref{a2}) and the functions
$F(x,y)$, $F_1(x,y)$ and $F_2(x,y)$, in~(\ref{a0}), and used these
to provide constructive methods to derive the linearizing
S-transformations for the case ${S_1}={S_2}= 0$ and the case
${S_1} \ne 0$ but ${S_3}={S_4}= 0$. This section provides an
alternative  procedure for constructing the f\/irst integrals and
Sundman transformations for S-linearizable equations, which can be
applied to both of the cases.

The key idea in this paper is that instead of f\/inding additional
relationships between the functions~$\phi$ and~$\psi$, we f\/ind
additional relationship between the functions~$A$ and~$B$ in~(\ref{a4}) and the functions $F(x,y)$, $F_1(x,y)$ and $F_2(x,y)$, in
(\ref{a0}). Equation~(\ref{a4}) implies
\begin{gather*}
 {\left(\frac{ {B_y}-{A_x}}{A}\right)}_y =  {{F_1}_y - 2{F_2}_x},
 \end{gather*}
which leads to the following missing relationship
 \begin{gather}\label{a7}
 {B_y}-{A_x} =  A({F_1} - 2{h_x}),
 \end{gather}
 where
\begin{gather}\label{a8}
h = \int {{F_2}(x,y)dy}+g(x),
\end{gather}
and $g(x)$ can be determined using Theorem \ref{t3} in case that
the ODE is S-linearizable.

This missing equation jointly with (\ref{a3}) give a new compact
S-linearizability criterion  for ODE of the form (\ref{a0}), given
in Theorem \ref{t3}. The S-linearizability criterion is used to
provide explicit general solutions for the auxiliary functions $A$
and $B$ which can be directly utilized to obtain the f\/irst
integral of (\ref{a3}), given in Theorem \ref{t4}, and hence the
Sundman transformations can be constructed using Theorem \ref{t1}.
Thus an alternative  procedure for constructing the f\/irst integral
and Sundman transformation is obtained.
\begin{Theorem}\label{t3}
Let us consider an equation of the form \eqref{a0}  and let $h =
\int {{F_2}(x,y)dy}+g(x)$. Equation~\eqref{a0}  is S-linearizable
if and only if
\begin{gather}\label{a9}
{{F_{1}}_x} + {F_1}{h_x} - h_x^2 - {h_{xx}} - {F_{y}} - {F}{F_2}=0,
\end{gather}
for some auxiliary function~$g(x)$.
\end{Theorem}

\begin{proof}
 Using the new relationship (\ref{a7}) with (\ref{a3}) one can get
 the following equations
\begin{gather}\label{a10}
{A_x} = A{h_x},\qquad
{A_y}=A{F_2},\qquad
{B_x}=AF,\qquad
{B_y} =  A({F_1} - {h_x}).
\end{gather}
The compatibility of the system (\ref{a10}),  i.e.~${A_{xy}} =
{A_{yx}}$ and ${B_{xy}} = {B_{yx}}$, leads to the criteria~(\ref{a9}).

In order to show that the new criterion~(\ref{a9}) is equivalent
to the one given in Theorem~\ref{t2}, we note that the system
consisting of equation~(\ref{a9}) and the second derivatives of
$h$ given by~(\ref{a8})
\begin{gather}\label{a11}
{h_{xx}}={{F_{1}}_x} + {F_1}{h_x} - h_x^2  - {F_{y}} - {F}{F_2},\qquad
{h_{yx}}={F_2}_x,\qquad
{h_{yy}}={F_2}_y.
\end{gather}
is compatible, i.e.~$h_{xy} = h_{yx}$, $h_{xxy} = h_{yxx}$ and
$h_{yyx} = h_{yxy}$, when the following equation holds
\begin{gather}\label{a12}
{h_x}({F_{1y}} - 2{F_{2x}}) + {F_{2x}}{F_1} + {F_{1xy}} -{F_{2xx}} - {F_{yy}} - {F_{y}}{F_2} - {F}{F_{2y}} = 0.
\end{gather}
Now, using ${S_1}$ and ${S_2}$ def\/ined by (\ref{a5}), equation~(\ref{a12}) can be rewritten as
\begin{gather*}
{S_1}{h_x} = {S_2} + {F_1}{S_1}.
\end{gather*}
Then clearly  if  ${S_1}=0$, then ${S_2}=0$ and if ${S_1} \ne 0$,
then
\begin{gather}\label{a14}
 {h_x} = \frac{{{S_2}}}{{{S_1}}} + {F_1}.
\end{gather}
Finally, substituting (\ref{a14}) in (\ref{a11}), gives
\begin{gather*}
{S_3}(x,y) = {\left( {\frac{{{S_2}}}{{{S_1}}}} \right)_y} -({F_{2x}} - {F_{1y}})=0,\\
{S_4}(x,y) = {\left( {\frac{{{S_2}}}{{{S_1}}}} \right)_x} +{\left( {\frac{{{S_2}}}{{{S_1}}}} \right)^2} + {F_1}\left({\frac{{{S_2}}}{{{S_1}}}} \right) + {F}{F_2} + {F_{y}}=0. \tag*{\qed}
\end{gather*}
\renewcommand{\qed}{}
\end{proof}
\begin{Remark}\label{r0}
In the case $S_1=S_2=0$, the criteria (\ref{a9}) can be
transformed by the change of variable
\begin{gather}\label{alt1}
g(x)=\ln{q(x)}+k(x),
\end{gather}
where $k'(x)=\frac{1}{2}F_1-\int {F_{2x}dy}$ to the well-def\/ined
ODE, equation~(\ref{eq-21}), in Theorem~\ref{th-1} and
$P=h-\ln{q}$ verif\/ies the system~(\ref{eq-20}) in Theorem~\ref{th-1}.

Moreover, in the case $S_3=S_4=0$, the criteria~(\ref{a9}) implies
equation~(\ref{a14}) which shows that $P=h$ verif\/ies the system~(\ref{eq-25}) in Theorem \ref{th-1} and hence solving this system
provides a~well-def\/ined ODE
\begin{gather}\label{alt2}
g'(x)=k(x),
\end{gather}
where $k(x)=\frac{{{S_2}}}{{{S_1}}} + {F_1}-\int {F_{2x}dy}$.

Hence when equation (\ref{a0}) is S-linearizable, one can solve
the criteria (\ref{a9}) for a function $g(x)$ by considering both
of~$x$ and~$y$ as independent variables. Or equivalently one can
get the function~$g(x)$ by equation~(\ref{alt1}) when $S_1=S_2=0$
whereas when $S_3=S_4=0$, the function $g(x)$ can be obtained from
equation~(\ref{alt2}).
\end{Remark}

In the next theorem, the general solution of the f\/irst integral is
given explicitly in terms of the function $h(x)$ where $h = \int
{{F_2}(x,y)dy}+g(x)$. It can be verif\/ied that this solution
coincides with the solution of the systems given in Theorem~\ref{th-1}. Hence, it provides an alternate direct procedure for
constructing the f\/irst integrals and the S-transformations.
\begin{Theorem} \label{t4}
Let us assume that equation \eqref{a0} is S-linearizable. Then
\eqref{a0}  has the first integral $w(x,y,y')=A(x,y)y'+B(x,y)$
where $A$ and $B$ are given by
\begin{gather*}
A(x,y) =  {{\rm{e}}^h}, \qquad
B(x,y)  =\int F   {{\rm e} ^{h  }}{dx}+\int  \left( {{\rm e}^{h}}(F_1 - h_x ) -\int  {{\rm e}^{h  }}( F_y  +F   F_2   ) {dx} \right){dy}, \end{gather*}
and $h$ is given by \eqref{a8}.
\end{Theorem}

\begin{proof}
The functions  $A$ and $B$ def\/ined by (\ref{a4}) can be
given explicitly by f\/inding the general solution  of the system~(\ref{a10}) where the second and the third equations of the system~(\ref{a10}) have general solution{\samepage
\begin{gather}\label{a18}
A =v  ( x  ) {{\rm e}^h},\qquad
B =\int F  A  {dx}+z ( y  ),
\end{gather}
for arbitrary functions $v  ( x  )$ and $z ( y
 )$.}

Substituting (\ref{a18}) in the f\/irst and the fourth equations of
the system (\ref{a10}) gives
\begin{gather}\label{a19}
v  ( x  ) = c_1,\qquad
z(y) =  \int \left(A({F_1} - {h_x})-\int (F  A_y+F_y A)  {dx}\right){dy} + k(x),
\end{gather}
for arbitrary functions $k ( x  )$.

Now, dif\/ferentiating (\ref{a19}) with respect to $x$  and using
the criterion (\ref{a9}) gives
\begin{gather*}
k_x =-\int A\big({{F_{1}}_x} + {F_1}{h_x} - h_x^2 - {h_{xx}} - {F_{y}} - {F}{F_2} \big){dy}=0 .
\end{gather*}
So $k(x)=c_2$, and f\/inally, from Theorem \ref{t1}, equation
(\ref{a0}) has the f\/irst integral $w(x,y,y')=A(x,y)y'+B(x,y)$ and
without loss of generality, one can choose $c_1=1$ and $c_2=0$ by
relabeling of $\frac{w(x,y,y')-c_2}{c_1}$.
\end{proof}

\begin{Remark}\label{r1}
An algorithmic implementation of our approach can be carried out
as summarized below.
Given a S-linearizable ODE of the form~(\ref{a0}). Use Theorem
\ref{t3} to determine an auxiliary function $g(x)$. Find the f\/irst
integrals using $g(x)$ and Theorem~\ref{t4}. Construct the Sundman
transformations~(\ref{OP}) using the f\/irst integral and Theorem~\ref{t1}. Since the free particle equation $u_{tt}=0$ has the
general solution $u(t)=c_1+c_2t$, f\/inally using the Sundman
transformations leads to the second integral of the ODE~(\ref{a0})
\begin{gather*}
\psi(x,y)=c_1+c_2 \mu(x),
\end{gather*}
where $t=\mu(x)$ is a solution of the f\/irst-order ODE
\begin{gather*}
\frac{dt}{dx}=\phi(x,\gamma(x,t)),
\end{gather*}
and  $y=\gamma(x,t)$ can be obtained by solving $c_1+c_2t =\psi
 ( x,y )$ for $y$.

But in case  that $\phi=\phi(x)$, $\mu(x)=\int{\phi(x)}dx$ and so
the Sundman transformation is a point transformation and leads to
the general solution~\cite{Warisa}.
\end{Remark}

In the next two examples,  we apply our approach to construct the
f\/irst integrals and use these to recover the Sundman
transformations of Muriel and Romero in~\cite{MR2010a}. In
addition, we provide the two-parameter family of solution in the
f\/irst example.
\begin{Example} \label{ex1}
Consider the ODE for the variable frequency oscillator~\cite{Mimura}
\begin{gather}\label{a27}
y'' + y{y'^2} = 0,
\end{gather}
Theorem \ref{t2} shows that the coef\/f\/icients of the equation
satisfy $S_1=0, S_2=0$. By Theorem~\ref{t3}, ODE (\ref{a27}) is
S-linearizable if and only if
\begin{gather}\label{a28}
g'' +{g'^2} = 0
\end{gather}
for some auxiliary function  $g(x)$. A particular solution of~(\ref{a28}) is $g(x)= \ln{x}$  so by (\ref{a8}) we have $h =
\frac{{{y^2}}}{2}+\ln{x}$  and hence using Theorem~\ref{t4}, we
can get the f\/irst integral $w(x,y,y')=A(x,y)y'+B(x,y)$ where
\begin{gather*}
A(x,y) = x\exp\left(\frac{y^2}{2}\right), \qquad
B(x,y)  =-\int \exp\left(\frac{y^2}{2}\right) dy.
\end{gather*}
Finally, the Sundman transformations can be constructed using
Theorem~\ref{t1} as follows{\samepage
\begin{gather*}
\psi(x,y) =\eta(I(x,y)) ,\qquad
\phi(x,y) =\frac{1}{x^2}\eta'(I(x,y)),
\end{gather*}
where $I(x,y)=x^{-1}  \int \exp\left(\frac{y^2}{2}\right) dy $.}

Now choosing $\eta(I)=I$ makes $\phi(x,y)=\phi(x)$ and so using
Remark~\ref{r1}  gives the two-para\-me\-ter  family of solutions of the
ODE~(\ref{a27})
\begin{gather*}
\text{erf\/i}\left(\frac{y}{\sqrt{2}}\right)=C_1 x+C_2,
\end{gather*}
where ${\rm erf\/i}(y)=\frac{2}{\sqrt{\pi}} \int_{0}^{y}{{\rm e}^{t^2}}dt$ is
the imaginary error function.
\end{Example}

\begin{Example}
Consider the equation
\begin{gather}\label{a213}
y'' - \left( {\tan {y} + \frac{1}{y}} \right){y'^2} + \left(
{\frac{1}{x} - \frac{{\tan {y}}}{{xy}}} \right)y' - \frac{{\tan
{y}}}{{{x^2}}} = 0.
\end{gather}
Theorem \ref{t2} shows that the coef\/f\/icients of this equation
satisfy ${S_1} \ne 0$ but ${S_3} = {S_4} = 0$. In addition, it
follows from Theorem \ref{t3} that ODE (\ref{a213}) is
S-linearizable  if and only if
\begin{gather}\label{a283}
g'' +{g'^2}-\left(\frac{1}{x}-\frac{\tan y}{xy}\right)g' = 0
\end{gather}
for some auxiliary function $g(x)$. The only solution of~(\ref{a283}) is $g(x) = C$ so by (\ref{a8}) we have $ h =
\ln{\left(\frac{\cos y}{y}\right)}$ and hence using Theorem~\ref{t4}, we can get the f\/irst integral
$w(x,y,y')=A(x,y)y'+B(x,y)$ where
\begin{gather*}
A(x,y) = \frac{{\cos y}}{y}, \qquad
B(x,y) = \frac{{\sin y}}{{xy}}.
\end{gather*}
Finally, the Sundman transformations can be constructed using
Theorem~\ref{t1} as follows
\begin{gather*}
\psi(x,y) =\eta(I(x,y)) ,\qquad
\phi(x,y) =xy \eta'(I(x,y)),
\end{gather*}
where $I(x,y)=x\sin{y}$.

One can show that there is no  $\eta(I)$ which makes
$\phi=\phi(x)$ and so using Remark~\ref{r1} for $\eta(I)=I$ gives
the two-parameter family of solution of the ODE (\ref{a213})
\begin{gather*}
x\sin{y}=c_1 +c_2\mu(x),
\end{gather*}
where the function $t=\mu(x)$ is a solution of the equation
\begin{gather*}
\frac{dt}{dx}=x\sin^{-1}\left(\frac{c_1+c_2t}{x}\right).
\end{gather*}
For example, if $c_2=0$, then one obtains the solution of ODE~(\ref{a213}) as
\begin{gather*}
x\sin{y}=c_1.
\end{gather*}
\end{Example}

As another application, we solve geodesic equations for surfaces
of revolution of constant curvature in a unif\/ied manner. Consider
a surface of revolution with parameterization $ ( f(y) \cos{x}
$, $f(y) \sin{x}, g(y)  )$ obtained by revolving the unit speed
curve $ ( f(y), g(y)  )$. The geodesic equations are~\cite{Pressley}
\begin{gather*}
\ddot{y}=f(y)f'(y)\dot{x}^2,\qquad
\frac{d}{dt}\big({f(y)}^2\dot{x}\big)=0,
\end{gather*}
where $\dot{y}=\frac{dy}{dt}$ and $\dot{x}=\frac{dx}{dt}$.

Using the formulas $\frac{dy}{dx}=\frac{\dot{y}}{\dot{x}}$ and
$\frac{d^2y}{dx^2}=\frac{\dot{x}\ddot{y}-\dot{y}\ddot{x}}{\dot{x}^3}$
gives
\begin{gather}\label{a21}
y'' -2\frac{f'(y)}{f(y)}{y'^2}-f'(y)f(y) = 0,
\end{gather}
which as special case for $f(y)=\sin{y}$ includes the equation for
geodesics on unit sphere given in \cite{Bluman2002,Stephani1989}.
\begin{Example} \rm \label{ex2}
We consider the nonlinear second-order ODE (\ref{a21}) for
$f(y)=y$, $f(y)=b+y$, $f(y)=\sin{y}$ and $f(y)=\sinh{y}$ describing the
geodesics on cone, plane, sphere and surface of conic type
respectively.

For the surfaces under consideration we have
$f'^2(y)-f(y)f''(y)=1$. It can be checked from Theorem~\ref{t2}
that the coef\/f\/icients of the equation satisfy $S_1=0$, $S_2=0$. In
addition, it follows from Theorem \ref{t3} that ODE~(\ref{a21})
for each of $f(y)=y$, $f(y)=b+y$, $f(y)=\sin{y}$ and $f(y)=\sinh{y}$ is
S-linearizable  if and only if
\begin{gather}\label{a22}
g'' +{g'^2} + 1=0,
\end{gather}
for some auxiliary function $g(x)$. A particular solution of
(\ref{a22}) is $g(x)= \ln{(\sin x)}$, so by (\ref{a8}) we have $ h
= \ln{\left(\frac{\sin x}{f^2(y)}\right)}$ and hence using Theorem~\ref{t4}, we can get the f\/irst integral
$w(x,y,y')=A(x,y)y'+B(x,y)$ where
\begin{gather*}
A(x,y) = \frac{\sin x}{f^2(y)}, \qquad
B(x,y)  =\frac{f'(y)}{f(y)} \cos x.
\end{gather*}
Finally, the Sundman transformations can be constructed using
Theorem~\ref{t1} as follows
\begin{gather*}
\psi(x,y) =\eta(I(x,y)) ,\qquad
\phi(x,y) =\left(\frac{f(y)}{f'(y)}\right)^2  \eta'(I(x,y)),
\end{gather*}
where  $I(x,y)=\frac{f(y)}{f'(y)}  \sin x$.

 Now choosing
$\eta(I)=\frac{1}{I}$ makes $\phi(x,y)=\phi(x)$ and so using
Remark \ref{r1}  gives the two-pa\-ra\-me\-ter  family of solution of the
ODE (\ref{a21})
\begin{gather*}
c_1 f(y)\sin x +c_2f(y)\cos x=f'(y).
\end{gather*}
\end{Example}

\begin{Example}  \label{ex3}
We consider the nonlinear second-order ODE
\begin{gather}\label{a35}
y'' -2\tanh{y} {y'^2}-\cosh{y}\sinh{y} = 0,
\end{gather}
that describes the geodesics on hyperboloid of one sheet.

Theorem~\ref{t2} shows that the coef\/f\/icients of the equation
satisfy $S_1=0$, $S_2=0$. It follows from Theorem~\ref{t3} that ODE~(\ref{a35}) is S-linearizable  if and only if
\begin{gather}\label{a36}
g'' +{g'^2} - 1=0,
\end{gather}
for some auxiliary function $g(x)$. A particular solution of~(\ref{a36}) is $g(x)= \ln{(\sinh x)}$, so by~(\ref{a8}) we have $
h = \ln{\left(\frac{\sinh x}{\cosh^2{y}}\right)}$ and hence using
Theorem~\ref{t4}, we can get the f\/irst integral
$w(x,y,y')=A(x,y)y'+B(x,y)$ where
\begin{gather*}
A(x,y) = \frac{\sinh x}{\cosh^2{y}}, \qquad
B(x,y)  =-\frac{\sinh{y}}{\cosh{y}} \cosh x.
\end{gather*}
Finally, the Sundman transformations can be constructed using
Theorem \ref{t1} as follows
\begin{gather*}
\psi(x,y) =\eta(I(x,y)) ,\qquad
\phi(x,y) ={\rm csch}^2 x  \eta'(I(x,y)),
\end{gather*}
where  $I(x,y)=\frac{\tanh{y}}{\sinh{x}}$.

 Now choosing
$\eta(I)=I$ makes $\phi(x,y)=\phi(x)$ and so using Remark~\ref{r1}
gives the two-para\-me\-ter  family of solution of the ODE~(\ref{a35})
\begin{gather*}
c_1\cosh{y}\sinh{x}-c_2\cosh{y}\cosh  x=\sinh{y}.
\end{gather*}
\end{Example}

\begin{Example}   \label{ex4}
We consider the nonlinear second-order ODE
\begin{gather}\label{a40}
y'' -2{y'^2}-{\rm e}^{2y} = 0,
\end{gather}
that describes the geodesics on pseudosphere.

Theorem~\ref{t2} shows that the coef\/f\/icients of the equation
satisfy $S_1=0$, $S_2=0$. It follows from Theorem~\ref{t3} that ODE~(\ref{a40}) is S-linearizable  if and only if
\begin{gather}\label{a41}
g'' +{g'^2}=0,
\end{gather}
for some auxiliary function $g(x)$. A particular solution of~(\ref{a41}) is $g(x)= \ln{x}$, so by~(\ref{a8}) we have $ h =
\ln{x}-2y$ and hence using Theorem~\ref{t4}, we can get the f\/irst
integral $w(x,y,y')=A(x,y)y'+B(x,y)$ where
\begin{gather*}
A(x,y) = x {\rm e}^{-2y},  \qquad
B(x,y)  =\frac{1}{2}\big({\rm e}^{-2y}-x^2\big).
\end{gather*}
Finally, the Sundman transformations can be constructed using
Theorem~\ref{t1} as follows
\begin{gather*}
\psi(x,y) =\eta(I(x,y)) ,\qquad
\phi(x,y) =-\frac{2}{x^2}  \eta'(I(x,y)),
\end{gather*}
where  $I(x,y)=\frac{{\rm e}^{-2y}}{x}+x$.

 Now choosing
$\eta(I)=I$ makes $\phi(x,y)=\phi(x)$ and so using Remark~\ref{r1}
gives the two-para\-me\-ter family of solutions of the ODE~(\ref{a40})
\begin{gather*}
{\rm e}^{-2y}+x^2=c_1x+2c_2.
\end{gather*}
\end{Example}

\section{Conclusion}
The recent Muriel--Romero characterization, Theorem \ref{t1}, of
the class of S-linearizable equations identif\/ies these as the
class of equations that admit f\/irst integrals of the form
$A(x,y)y'+B(x,y)$.  In this paper, a new characterization of
S-linearizable equations in terms of the coef\/f\/icients and one
auxiliary function is given in Theorem~\ref{t3}. This criterion is
used to directly provide explicit general solutions for the
auxiliary functions~$A$ and~$B$ Theorem~\ref{t4}. So, using
Theorem~\ref{t1}, the linearizing generalized Sundman
transformations can be constructed by solving the f\/irst-order ODE~(\ref{a42}). Finally, it is shown in~\cite{MR2010a} that an
equation of the form~(\ref{a0}) is S-linearizable and linearizable
via a point transformation if and only if $S_1=S_2=0$. It is also
known that the generalized Sundman transformation is a point
transformation if and only if $\phi=\phi(x)$. So, by Remark~\ref{r1}, the generalized Sundman transformation leads to the
general solution $\psi(x,y)=c_1+c_2  \int{\phi(x)}dx$  if and only
if $S_1=S_2=0$.

Our method is illustrated in examples where we recover the Sundman
transformations of Muriel and Romero in~\cite{MR2010a}.
Furthermore, the system of geodesic equations for surfaces of
revolution is expressed  as a single second-order ODE. It is
noticed that this ODE is S-linearizable for surfaces of revolution
with constant curvature. The method is applied to f\/ind the general
solution of these geodesics  in a unif\/ied manner.

\subsection*{Acknowledgments}

  The authors would like to thank the King Fahd University of
Petroleum and Minerals for its support and excellent research
facilities. They also thank the reviewers for their comments which
have considerably improved the paper.

\pdfbookmark[1]{References}{ref}
\LastPageEnding

\end{document}